\documentclass[11pt, a4paper, headings=standardclasses, openany, bibliography=totocnumbered]{amsart}
\usepackage{standalone}
\usepackage[T1,T2A]{fontenc}
\DeclareUnicodeCharacter{0306}{HERE!HERE!}
\usepackage[automark, autooneside=false]{scrlayer-scrpage}  
\usepackage[usenames,dvipsnames]{xcolor}
\usepackage[T1]{fontenc}
\usepackage[utf8]{inputenc}
\usepackage{lmodern}
\usepackage[UKenglish]{babel}
\usepackage{amsmath, amssymb, amsfonts, amsthm, mathtools, bbm}
\usepackage{extarrows}
\usepackage{mathrsfs}
\usepackage{pgfplots}
\usetikzlibrary{arrows}
\usepackage{tikz-cd}
\usepackage{graphicx}
\usepackage{faktor}
\usepackage{young}
\usepackage[centertableaux]{ytableau}
\usepackage{verbatim} 
\usepackage{comment}
\usepackage{wasysym}
\usepackage{enumitem}
\usepackage{microtype}

\usepackage{csquotes}
\usepackage{hyperref}
\usepackage{cleveref}	
\usepackage[style=alphabetic,backend=biber]{biblatex}
\usepackage[normalem]{ulem}
\usepackage{color}      
\usepackage{xparse}
\usepackage{longtable}

\usepackage{braids}
\usepackage{adjustbox}
\usepackage{forest} 
\usepackage{dynkin-diagrams}
\usepackage{rank-2-roots}
\definecolor{light-gray}{gray}{0.85}

\binoppenalty=\maxdimen         
\relpenalty=\maxdimen           

\newcommand{\opn}[1]{\operatorname{#1}}

\newtheorem{theo}{Theorem}
\numberwithin{theo}{section}
\newtheorem{lemm}[theo]{Lemma}

\newtheorem{prop}[theo]{Proposition}
\newtheorem{koro}[theo]{Corollary}

\theoremstyle{definition}
\newtheorem{defi}[theo]{Definition}

\newtheorem{beis}[theo]{Example}
\newtheorem{rema}[theo]{Remark}

\newtheorem{nota}[theo]{Notation}

\newcommand{\odel}{\bar{\del}} 
\newcommand{\botimes}{\textcolor{blue}{\otimes}} 
\newcommand{\cotimes}{\textcolor{Turquoise}{\otimes}} 
\newcommand{\rotimes}{\textcolor{red}{\boxtimes}} 
\DeclareMathOperator{\red}{red} 
\DeclareMathOperator{\Hom}{Hom}
\newcommand{\mix}{\opn{mix}}

\newcommand{\esi}{\prescript{}{s}{\underline{\hat{i}}}} 
\newcommand{\esm}{\prescript{}{s}{\underline{\hat{m}_{st}}}} 
\newcommand{\eti}{\prescript{}{t}{\underline{\hat{i}}}} 
\newcommand{\etm}{\prescript{}{t}{\underline{\hat{m}_{st}}}} 
\newcommand{\sm}{\prescript{}{s}{\hat{m}}} 
\newcommand{\tm}{\prescript{}{t}{\hat{m}}} 
\DeclareMathOperator{\rR}{Rel} 
\DeclareMathOperator{\Der}{Der} 
\DeclareMathOperator{\Diff}{Diff} 
\DeclareMathOperator{\Ann}{Ann} 
\DeclareMathOperator{\Ext}{Ext}
\DeclareMathOperator{\End}{End}

\newcommand{\del}{{\partial}}
\renewcommand{\C}{\mathbb{C}}

\DeclareMathOperator{\triv}{triv}
\DeclareMathOperator{\id}{id}

\renewcommand{\U}{\mathbf{U}}

\newcommand{\Mod}{\text{-}\mathbf{Mod}}
\newcommand{\Bim}{\text{-}\mathbf{Bim}\text{-}}
\DeclareMathOperator{\Vect}{\mathbf{Vect}}

\newcommand{\ts}{\textsuperscript}

\DeclareMathOperator{\sl2}{\mathfrak{sl}_2}

\newcommand{\sln}{\mathfrak{sl}_n}

\newcommand{\h}{\mathfrak{h}}

\DeclareMathOperator{\gl}{\mathfrak{gl}}

\DeclareMathOperator{\M}{M}

\DeclareMathOperator{\Sym}{Sym}
\DeclareMathOperator{\Quot}{Quot}

\DeclareMathOperator{\op}{op}

\DeclareMathOperator{\catC}{\mathscr{C}}

\DeclareMathOperator{\nh}{n\mathcal{H}} 
\DeclareMathOperator{\nc}{n\mathcal{C}} 

\newcommand{\W}{\mathcal{W}} 

\DeclareMathOperator{\mU}{\mathbbm{1}}

\DeclareMathOperator{\mR}{\mathbb{R}}

\newcommand{\Hecke}{\mathcal{H}}

\newcommand{\shift}[1]{\langle#1\rangle} 

\newcommand\restr[2]{{
		\left.\kern-\nulldelimiterspace 
		#1 
		\vphantom{\big|} 
		\right|_{#2} 
}} 

\newcommand{\wo}{w_o} 
\newcommand{\blank}{{-}}

\DeclareMathOperator{\Gal}{Gal}

\DeclareFieldFormat[article,unpublished]{citetitle}{\mkitalic{#1}\isdot}
\DeclareFieldFormat[article,unpublished]{title}{\mkbibitalic{#1}\isdot}
\DeclareFieldFormat[article,book]{number}{\mkbibbold{#1}}
\DeclareFieldFormat[article]{volume}{\mkbibbold{#1}}

\usetikzlibrary{calc}
\usetikzlibrary{tqft}
\newcommand{\cbox}[1]{ \vcenter{ \hbox{#1} } }
\usetikzlibrary{patterns,decorations.pathreplacing}
\tikzset{
	tldiagram/.style={thick, scale=0.35}
}
\newcommand{\tlcoord}[2]
{
	(2*#2 , 3*#1)
}
\newcommand{\makecdots}[2]
{
	\node at (2*#2, 3*#1 + 1.5) {$\cdots$};
}

\newcommand{\lineup}{-- ++(0,3)}
\newcommand{\halflineup}{-- ++(0,1.5)}

\newcommand{\linewave}[2]{
	.. controls +(0,1.5*#1) and +(0,1.5*-#1) .. ++(2*#2, 3*#1)
}

\newcommand{\timesr}{\times_{\! R}}
\newcommand{\rtimesr}{\textcolor{red}{\timesr}} 
\newcommand{\btimesr}{\textcolor{blue}{\timesr}}
\newcommand{\otimesr}{\otimes_{\! R}}

\newcommand{\dlineup}{\halflineup \onedot \halflineup}

\newcommand{\onedot}{node {$\bullet$}}

\newcommand{\dcaprightsmall}{arc (180:90:1) node {\scalebox{0.7}{\textbullet}} arc (90:0:1)}

\newcommand{\dcuprightsmall}{arc (180:270:1) node {\scalebox{0.7}{\textbullet}} arc (270:360:1)}

\newcommand{\negtwista}[2]{
	\draw [ultra thick] \tlcoord{0+3*#1}{0+2*#2} -- \tlcoord{0.45+3*#1}{0.4+2*#2};
	\draw [ultra thick] \tlcoord{0.55+3*#1}{0.6+2*#2} -- \tlcoord{1+3*#1}{1+2*#2};
	\draw  \tlcoord{0+3*#1}{1+2*#2} -- \tlcoord{1+3*#1}{0+2*#2};
}

\newcommand{\negtwist}[2]{
	\negtwista{#1}{#2}
	\negtwist{1+#1}{#2}
}

\newcommand{\predtlcupcap} 
{
	\begin{tikzpicture}[scale=0.092]
		\draw \tlcoord{1}{0} \dcuprightsmall;
		\draw \tlcoord{0}{0} \dcaprightsmall;
	\end{tikzpicture}
}

\setlist[enumerate]{leftmargin=*, label= \roman*)}

\setlist{
	listparindent=\parindent,
	parsep=0pt
}

\begin{document}
	\title[Nil Hecke algebras - an algebroidists perspective]{{A study of nil Hecke algebras via Hopf algebroids}}
	\author[Z. Wojciechowski]{Zbigniew Wojciechowski}
	\address{Z. W.: Institut für Geometrie, Technische Universität Dresden, Germany}
	\email{zbigniew.wojciechowski@tu-dresden.de}
	
	\maketitle
	\begin{abstract}
		Hopf algebroids are generalizations of Hopf algebras to less commutative settings. We show how the comultiplication defined by Kostant and Kumar turns the affine nil Hecke algebra associated to a Coxeter system into a Hopf algebroid without an antipode. The proof relies on mixed dihedral braid relations between Demazure operators and simple reflections. For researchers new to Hopf algebroids we include additional examples from ring theory, representation theory, and algebraic geometry. 
	\end{abstract}
\section{Introduction}
To understand an abstractly defined algebra $A$ over a field $k$, finding a faithful representation of $A$ on some vector space is helpful. If we think one categorical level higher, an algebra becomes a monoidal additive $k$-linear category $(\catC,\otimes, \mU)$ and a faithful representation becomes a faithful $k$-linear (strong) monoidal functor $F$ to some category we understand, examples being
\[
\begin{tikzcd}[column sep=small]
	& (\Vect_k,\otimes_{k},k)      & k \text{ a field,}                   \\
	F\colon (\catC,\otimes, \mU) \arrow[end anchor={[xshift=-1.25em]}]{ru} \arrow[r] \arrow[end anchor={[xshift=-2em]}]{rd} & (R\Mod, \otimesr, R)     & R \text{ a commutative $k$-algebra,} \\
	& (R\Bim R, \otimesr, R) & R \text { any $k$-algebra.}          
\end{tikzcd}
\]
We are interested in the case where $\catC=H\Mod$ for a $k$-algebra $H$. By \cite[Theorem 5.1]{schauenburg98}, the existence of $F$ together with the choice of target category corresponds to additional algebraic structure on $H$:
\[
\setlength{\arraycolsep}{11.5pt}
	\begin{array}{ccc}
		\Vect_k & H \text{ is a $k$-bialgebra,} \\
		R\Mod & H \text{ is a Sweedler $R$-bialgebroid, \cite[Definition 5.6]{sweedler74}} \\
		R\Bim R & H \text{ is a Takeuchi $R$-bialgebroid, \cite[Definition 4.5]{takeuchi77}} \\ 
	\end{array}
\]
We aim to find natural examples of both kinds of bialgebroids in representation theory. This paper focuses on Sweedler bialgebroids, which we call bialgebroids. Bialgebroids over $R$ are triples $(H,\Delta, \varepsilon)$, where $H$ is a $k$-algebra $H$ containing $R$ with a comultiplication $\Delta\colon H\to H \otimesr H$ and a counit $\varepsilon\colon H\to R$ satisying certain conditions (see \Cref{definition bialgebroid}), which generalize the notation of an $R$-bialgebra, without requiring $H$ to be an $R$-algebra. This is a generalization in three ways. First, the counit $\varepsilon$ is not required to be an algebra morphism. This allows for more interesting actions of $H$ on $\mU=R$. Second, the comultiplication $\Delta\colon H\to H \otimesr H$ is required to land in the Takeuchi product $H\timesr H\subseteq H\otimesr H$. This algebraic condition is often non-trivial to check and corresponds in natural examples to some weak commutativity between elements of $R$ and $H$. Finally, the existence of an antipode, which turns bialgebras into Hopf algebras, is not the condition one requires for bialgebroids to become Hopf algebroids. Instead, one requires the existence of an inverse for the Galois map. This inverse may be induced by one (or multiple!) antipodes $S\colon H\to H^{\op}$ or no antipode. 

\subsection*{Overview and results}
 \Cref{section on bialgebroid basics} recalls the definitions and some elementary facts on Hopf algebroids. Experts on these topics should come there for the notation, specifically \Cref{notation on colorful tensor products} is essential, where we explain what red and blue tensor symbols mean. 
 \Cref{section examples of bialgebroids} then contains some of the most important examples of bialgebroids coming from ring theory, representation theory, and algebraic geometry. We discuss how the definition of bialgebra relates to that of a bialgebroid (\Cref{example bialgebroid 1}), matrix algebras and semisimple algebras (Examples \ref{example bialgebroid 2}, \ref{example bialgebroid 3}), path algebras of quivers (\Cref{example bialgebroid 4}), and finally Weyl algebras and differential operators (Examples \ref{example bialgebroid 5},\ref{example bialgebroid 6}). Lastly, we mention two general constructions of bialgebroids --- endomorphism algebras under freeness assumptions in \Cref{example bialgebroids in finite free endo case}, and twisted group algebras in \Cref{example twisted group algebras}. Examples \ref{example bialgebroid 5}, \ref{example bialgebroids in finite free endo case}, \ref{example twisted group algebras} all provide helpful intuition for \Cref{section on Nil--Hecke basics}, which is on the new results. 
 
 In \Cref{section on Nil--Hecke basics} we show that the affine nil Hecke algebra introduced in \cite{kostant86} associated to a Coxeter system $(W,S)$ and its geometric representation $\h^{*}$ becomes a cocommutative Hopf algebroid over $R=\Sym(\h^{*})$ without antipode (\Cref{theorem nil hecke is hopf algebroid}, \Cref{theorem no antipode}). This structure relies on mixed dihedral braid relations between Demazure operators and simple reflections; see \Cref{theorem mixed relations}. We prove these relations indirectly via the embedding $\nh\hookrightarrow Q\star W$, where $Q$ is the fraction field of $R$, and using the comultiplication on $Q\star W$ from \cite{kostant86}.
 In the case where $W$ is finite, the Hopf algebroid structure induces the monoidal structure $\otimes_R$ on $\nh\Mod$, which makes the Morita equivalence $\nh\Mod \to \Sym(\h^{*})^{W}\Mod$ monoidal, see \Cref{corollary equivalence is monoidal}. Our methods 

\subsection*{Acknowledgments}
The author would like to thank the true algebroidists Myriam Mahaman and Ulrich Krähmer for discussions and guiding me through the Hopf algebroid literature. Special thanks go to Christian Lomp for asking me about the definition of the nil Hecke algebra. At that moment, I had the idea for the project set.
\section{Preliminaries on bialgebroids and Hopf algebroids} \label{section on bialgebroid basics}
This section is a recollection of definitions and facts from \cite{sweedler74}, written for representation theorists interested in using bialgebroids in their research. We start with fixing notations.
\begin{nota}
	All rings/algebras are associative and unital, in general not commutative. Ring and algebra morphisms are unital. Throughout this paper $k$ denotes a field, $R$ a commutative $k$-algebra, $H$ a $k$-algebra with an embedding of algebras $R\subseteq H$. We regard $H$ as $R$-$R$-bimodule with the obvious $R$-bimodule structure $r\cdot h \cdot r'=rhr'$ where $r,r'\in R, h\in H$. We set $H^{\op}\coloneqq H$ both as a $k$-vector space and $R$-$R$-bimodule, and view it as $k$-algebra via the opposite product $\cdot_{\op}$. In general, the bimodule structure and opposite product are not compatible; for instance, for $r\in R, h,h'\in H^{\op}$ we have $r\cdot (h\cdot_{\op} h')=rh'h$, while $(r\cdot h)\cdot_{\op} h'=h'rh$.
\end{nota} 
\begin{nota}[Blue and red tensor products] \label{notation on colorful tensor products}
	Given $M=\prescript{}{R}{M}, N=\prescript{}{R}{N}$ two left $R$-modules we will write
	$M\botimes N\coloneqq \prescript{}{R}{M}\otimesr\prescript{}{R}{N}$ for the tensor product of left $R$-modules. If $M=M_R$ is a right $R$-module, we write
	$M\rotimes N\coloneqq M_{R}\otimesr\prescript{}{R}{N}$ for the balanced tensor product. The same color scheme applies to tensoring morphisms and elements. As a mnemonic compare $ \text{\textcolor{blue}{b\underline{l}ue  circ\underline{l}e}}\simeq \text{\underline{l}eft}$, $\text{\textcolor{red}{\underline{r}ed squa\underline{r}e}}\simeq \text{\underline{r}ight}$, 
	where the left/right refers to the module $M$.	
\end{nota}

\begin{defi} \label{definition takeuchi products}
	We define subsets $H\btimesr H\subseteq H\botimes H$ and $H\rtimesr H^{\op}\subseteq H\rotimes H^{\op}$ by setting 
	\begin{align*}
		H\btimesr H&\coloneqq\!\left\{ \sum_{i}h_i \otimes h_i'\mid \forall r\in R\colon \sum_{i}(h_ir) \otimes h_i'= \sum_{i}h_i \otimes (h_i'r)\right\}, \\
		H\rtimesr H^{\op}&\coloneqq\!\left\{ \sum_{i}h_i \otimes h_i'\mid \forall r\in R\colon \sum_{i}(rh_i) \otimes h_i'= \sum_{i}h_i \otimes (h_i'r)\right\}.
	\end{align*}
	Both are called Takeuchi products.
\end{defi}
The following lemma is straightforward to check.
\begin{lemm} \label{lemma takeuchi products}
	The following hold:
	\begin{enumerate}
		\item Both $H\btimesr H$ and $H\rtimesr H^{\op}$ become $k$-algebras via the component-wise multiplication, moreover they contain $R$ as subalgebra $R\cdot 1\botimes 1\subseteq H\botimes H$ respectively $R\cdot 1 \rotimes 1\subseteq H\rotimes H^{\op}$.
		\item Let $M, N$ be two left $H$-modules. Then $M\botimes N$ is canonically a left $H\btimesr H$-module and $\Hom_R(M,N)$ a left $H\rtimesr H^{\op}$-module.
	\end{enumerate}
\end{lemm}
The following remark motivates \Cref{definition takeuchi products} and \Cref{lemma takeuchi products}.

\begin{rema} \label{remark why ii is different}
	Let $M, N$ be two left $H$-modules. We want to turn $M\botimes N$ respectively $\Hom_R(M,N)$ into left $H$-modules. However, first, we want to turn them into modules over some tensor product $H\otimes H$ or $H\otimes H^{\op}$ of algebras and then, in the second step, pull back this action to $H$ via some algebra morphism. At first sight, the candidates to act on these spaces are $H\botimes H$ respectively $H\rotimes H^{\op}$. However, the naive component-wise multiplication is not well-defined for both. The Takeuchi products make this idea work.
\end{rema}

\begin{defi} \label{definition bialgebroid}
	We call $H$ a left bialgebroid if it is equipped with morphisms of left $R$-modules $\Delta\colon H\to H\botimes H, \varepsilon\colon H\to R$, such that:
	\begin{enumerate}[itemsep=1pt]
		\item The triple $(H, \Delta, \varepsilon)$ is a coassociative, counital coalgebra in left $R$-modules. 
		\item \label{item v Takeuchi product Hopf algebroid} The comultiplication maps into the Takeuchi product  $\Delta\colon H\to  H\btimesr H$ and is a $k$-algebra morphism when regarded in this way.
		\item \label{item iv action on R Hopf algebroid} The map 
		\begin{equation} \label{definition rho epsilon}
			\rho_{\varepsilon} \colon H\to \End_k(R), \quad  h\mapsto (r \mapsto h(r)\coloneqq\varepsilon(hr))
		\end{equation} turns $R$ into a left $H$-representation such that the restriction $\restr{\rho_\varepsilon}{R}$ is the regular representation $\rho_{\operatorname{reg}}\colon R\to \End_k(R)$ of $R$ on $R$. 
	\end{enumerate}
	A bialgebroid is called Hopf or simply Hopf algebroid if the Galois map given by the composition 
	\[
	\Gal_H\colon H\rotimes H\xrightarrow{\Delta\,\rotimes \id_H}  H\botimes H\rotimes H \xrightarrow{\id_H \botimes \mu} H\botimes H
	\]
	is bijective, where $\mu\colon H\rotimes H\to H$ is the multiplication map. We define the red map $\red_H\colon H\to H\rotimes H$ by $\red_H(h)\coloneqq\Gal_H^{-1}(h\otimes 1)$. 
\end{defi}
\begin{rema}
	Condition \ref{item v Takeuchi product Hopf algebroid} is a natural assumption on $\Delta$ by the discussion in \Cref{remark why ii is different}. Condition \ref{item iv action on R Hopf algebroid} is a weakening of $\varepsilon\colon H\to R$ being an $R$-algebra morphism. In ring theoretic terms, \ref{item iv action on R Hopf algebroid} can be rephrased as $\ker(\varepsilon)\subseteq H$ being a left ideal and the counit fixing $R$, as formula $\restr{\varepsilon}{R}=\id_R$.
\end{rema}
\begin{rema}
	The map $\red_H$ associated with a Hopf algebroid $H$ is called translation map in the literature. One can check that $\red_H$ is an algebra morphism $\red_H\colon H\to H\rtimesr H^{\op}$. Since $\Gal_H$ is right $H$-linear, the inverse $\Gal_H^{-1}$ is also right $H$-linear, if it exists. Hence $\Gal_H^{-1}$ is uniquely determined by $\red_H$, which can be defined on generators.
\end{rema}
\begin{lemm} \label{lemma module category becomes monoidal} 
	Let $(H, \Delta, \varepsilon)$ be an $R$-bialgebroid. The category $H\Mod$ of left $H$-modules becomes monoidal via the tensor product $\botimes$ with unit $\mU=(R,\rho_{\varepsilon})$. The forgetful functor 
	$\opn{res}^{H}_R\colon H\Mod\to R\Mod$ is monoidal. If $H$ is Hopf the space $\Hom_R(M,N)$ becomes an $H$-module and an internal $\Hom$ in $(H\Mod, \botimes, \mU)$.
\end{lemm}
\begin{proof}
	See \cite[Theorem 3.5]{schauenburg00} for the original result and \cite[Remark 4.26]{kraehmermahaman24} for a short discussion of this fact using notation more similar to ours.
\end{proof}
\section{Examples of bialgebroids} \label{section examples of bialgebroids}
\begin{beis} \label{example bialgebroid 1}
	An $R$-bialgebra is exactly an $R$-bialgebroid which satisfies $R\subseteq \operatorname{Z}(H)$. In this case $H\btimesr H = H\botimes H$. In this sense, a bialgebroid is a less commutative object than a bialgebra.  Additionally a bialgebra is a Hopf algebra, i.e.\ admits an antipode $S\colon H\to H$, if and only if it as a bialgebroid is a Hopf algebroid. Indeed given a Hopf algebra one sets $\red_H(h)=h_{(1)}\rotimes \,S(h_{(2)})$, where $\Delta(h)=x_{(1)}\botimes \, x_{(2)}$ in Sweedler notation. On the other hand given $\red_H$ one can recover the antipode $S$ as
	\[
	S=(\varepsilon\rotimes \id_H)\circ \red_H.
	\]
	Note that this is only well-defined since $\varepsilon$ becomes right $R$-linear in this case.
\end{beis}
\begin{beis} \label{example bialgebroid 2}
	Consider the embedding $R=k^n\subseteq H=\M_{n\times n}(k)$ of diagonal matrices in all $n\times n$-matrices. Then $H$ is an $R$-bialgebroid via the comultiplication and counit 
	\begin{align*}
		\Delta\colon \M_{n\times n}(k)&\to \M_{n\times n}(k)\textcolor{blue}{\times_{k^n}} \M_{n\times n}(k), \quad E_{ij}\mapsto E_{ij}\botimes E_{ij} \\
		\varepsilon \colon \M_{n\times n}(k)&\to k^n, \quad E_{ij}\mapsto e_i\coloneqq E_{ii}
	\end{align*}
	where $(E_{ij})_{1\leq i,j\leq n}$ is the standard basis of $H$. The map $\rho_{\varepsilon}$ is the isomorphism $\M_{n\times n}(k)\cong \End_{k}(k^n)$. In this example, the choice of the algebra embedding of $k^n$ into $\M_{n\times n}(k)$ corresponds to a choice of basis of $k^n$, moreover it uniquely determines the comultiplication. The counit $\varepsilon$ evaluates a matrix at the identity matrix, viewed as the vector $(1,\ldots,1)^{\intercal}$. The map $\red_H\colon E_{ij}\mapsto E_{ij}\rotimes E_{ji}$ shows that $\M_{n\times n}(k)$ is even a Hopf algebroid. The monoidal structure on $\M_{n\times n}(k)$-modules is built in such a way that the Morita equivalence
	\begin{equation*}
		\Vect_k\to \M_{n\times n}(k)\Mod, \quad V\mapsto k^{n}\otimes_k V
	\end{equation*} becomes a closed monoidal equivalence.
\end{beis}
\begin{beis}\label{example bialgebroid 3}
	Assume $k=\bar{k}$. Similarly to \Cref{example bialgebroid 2}, every semisimple $k$-algebra $H$ is an $R$-Hopf algebroid, where $R\subseteq H$ corresponds to the product of diagonal matrix subalgebras under an Artin--Wedderburn isomorphism
	\[
	\prod_{i=1}^{r}k^{n_i}\cong R \subseteq H\cong \prod_{i=1}^r \M_{n_i\times n_i}(k),
	\]
	where $\Delta, \varepsilon$ can be calculated explicitly, provided one knows the isomorphism. In this case, $\Delta$ defines an isomorphism $H\to H\btimesr H$ of algebras. Tensoring a representation with an irreducible one over $R$ returns the corresponding isotypical component of the given representation. As a special case, the Iwahori--Hecke algebra $\Hecke_q(S_n)$ for generic $q\in \C^{\times}$ (that is $q$ is not a $(2m)$\ts{th} root of unity for all $0\leq m\leq n$) can be turned into a Hopf algebroid over its Gelfand--Zetlin subalgebra, that is the subalgebra generated by Jucys--Murphy elements. 
\end{beis}
\begin{beis} \label{example bialgebroid 4}
	For any quiver $Q=(Q_0,Q_1,s,t)$ the path algebra $H=kQ$ becomes a bialgebroid over $R=kQ_0$. The comultiplication $\Delta$ maps each path $\gamma$ to $\gamma\botimes \gamma$, and the counit maps each $\gamma$ to $e_{t(\gamma)}$ the idempotent constant path corresponding to the target vertex of $\gamma$. This bialgebroid structure corresponds on the representation theory of quivers side to the monoidal structure which tensors two quiver representations vertex-wise, that is for two quiver representations $V=((V_i)_{i\in Q_0}, (\varphi_{\alpha,V})_{\alpha\in Q_1}$ one has $(V\otimes W)_i= V_i\otimes W_i$ for $i\in Q_0$ and $\varphi_{\alpha, V\otimes W}=\varphi_{\alpha, V}\otimes \varphi_{\alpha, W}$. The monoidal unit $\mU$ is as vector space $kQ_0$, but not equipped with the trivial action, and instead equiped with the augmentation action $\rho_{\varepsilon}$. As quiver representation this is the one, where $k$ is assigned to each vertex and the identity is assigned to each arrow. Almost no path algebras of quivers are Hopf algebroids. If there is any arrow $\alpha\in Q_1$, then $\alpha\,\botimes\, e_{t(\alpha)}$ has no preimage under the Galois map. Hence, a path algebra $kQ$ is a Hopf algebroid if and only if $Q_1=\emptyset$, and hence $R=kQ_0=kQ=H$. 
\end{beis}
\begin{rema}
	We make two remarks on \Cref{example bialgebroid 4}. As discussed we have $\mU=kQ_0$ as a vector space, but $\mU$ is not the representation $kQ_0=\bigoplus_{i\in Q_0}S_i$, which appears when studying Koszul duality. This fact also becomes obvious by the Eckmann-Hilton argument, which shows that the $\Ext$-algebra of the monoidal unit $\mU$ becomes graded commutative, which in general does not apply to the Koszul dual $(kQ)^{!}$. Second, note that the bialgebroid structure on $kQ$ does not descend to a bialgebroid structure of quotients $kQ/I$. When $I$ consists of monomial relations, $\Delta$ is still well-defined, however $\varepsilon$ is not. Hence, $(kQ/I\Mod, \botimes)$ becomes a monoidal category without a unit, also called a semigroup category. 
\end{rema}
\begin{beis}  \label{example bialgebroid 5}
	Consider the Weyl algebra $\W_n$, that is the $k$-algebra with $2n$ generators $x_1,\ldots, x_n, \del_{1}, \ldots, \del_{n}$ and relations
	\[
	x_ix_j=x_jx_i, \quad \del_{i}\del_{j}=\del_{j}\del_{i}, \quad \del_{i}x_j=x_j\del_{i}+\delta_{i,j}. 
	\]
	Then $H=\W_n$ becomes a $k[x_1,\ldots, x_n]$-Hopf algebroid by setting
	\begin{gather*}
		\Delta\colon \W_n\to \W_n\textcolor{blue}{\times_{k[x_1,\ldots, x_n]} }\W_n, \quad \del_{i}\mapsto \del_i\botimes 1+ 1\botimes \del_i \\
		\qquad \varepsilon \colon \W_n\to k[x_1,\ldots, x_n],  x_{1}^{a_1}\cdots x_{n}^{a_n}\del_{1}^{b_1}\cdots \del_{n}^{b_n} \mapsto \delta_{0,b_1+b_2+\cdots +b_n}x_{1}^{a_1}\cdots x_{n}^{a_n}.
	\end{gather*}
	It is a standard but good exercise to check that $\Delta$ indeed maps into the Takeuchi-product. The corresponding representation $\rho_{\varepsilon}$ is the polynomial representation of $\W_n$, where the $\del_i$'s act by partial derivatives on $k[x_1,\ldots, x_n]$. We have $\red_{\W_n}(\del_i)=\del_i\rotimes1-1\rotimes \del_i$.
\end{beis}
\begin{beis} \label{example bialgebroid 6}	
	More generally, $H=\Diff(R)$, the algebra of differential operators on any smooth, complex, affine variety $X=\operatorname{Spec}(R)$, becomes an $R$-Hopf algebroid. Concretely $\Diff(R)$ is defined as 
	\[
	\Diff(R)\coloneqq \bigcup_{i\geq0}\Diff_i(R), \quad \Diff_i(R)\coloneqq\Ann(\ker(\mu)^{i+1})\subseteq \End_k(R)
	\]
	where $\mu\colon R\otimes_k R\to R$ is the multiplication map, and the annihilator is taken with respect to the canonical action of $R\otimes_k R= R\otimes_k R^{\op}$ on $\End_k(R)$. Under the assumption that $\operatorname{Spec}(R)$ is smooth, $\Diff(R)$ is generated by the multiplication operators $\Diff_0(R)\cong R$ and the derivations $\Der_k(R)\subseteq \Diff_1(R)$. Moreover the only additional relations are the Weyl relation
	\[
	r\del=\del r + \del(r), \quad \text{where } r\in R, \del\in \Der_k(R)
	\]	
	and the relation $(r\del)=r\del$ identifying the derivation $(r\del)$ with the product of $r\in R$ with $\del\in \Der(R)$ in $\Diff(R)$. The comultiplication is obtained by making every derivation of $R$ primitive. In different language $\Diff(R)$ is canonically isomorphic to $\U(R,\Der_k(R))$, the universal enveloping algebra of the Lie--Rinehart algebra $(R,\Der_k(R))$, see \cite[\S2, \S3]{kraehmermahaman24}
	for an overview. Other helpful references on this topic are \cite[\S 18]{sweedler74}, \cite{rinehart63}, and \cite{moerdijk10}. The non-smooth case is much more complicated, the results \cite{kraehmermahaman24} on singular curves suggest that $\Diff(R)$ still becomes a Hopf algebroid, however with a highly non-trivial comultiplication. 
\end{beis}
We end this section by discussing two other families of examples of bialgebroids, closely related to the nil Hecke case later on. The first is the running example in \cite{sweedler74}.
\begin{beis} \label{example bialgebroids in finite free endo case}
	Let $R'\subseteq R$ be commutative rings such that $R$ is finite free as $R'$-module. Famously, one has a Morita equivalence
	\[
	R\otimes_{R'}{-}\colon R'\Mod\overset{\simeq}{\longrightarrow} \End_{R'}(R)\Mod.
	\]
	Pulling the (symmetric) monoidal structure $\otimes_{R'}$ on $R'\Mod$ back through the equivalence defines a (symmetric) monoidal structure on $\End_{R'}(R)\Mod$, which agrees with $\botimes$ coming from the cocommutative $R$-bialgebroid structure on $\End_{R'}(R)$, which we describe in the following. First, we identify
	\[
	\End_{R'}(R)\cong R\otimes_{R'} R^{*}, \text{ where } R^{*}=\Hom_{R'}(R,R').
	\] 
	Since $R$ is a finite free algebra over $R'$, $R^{*}$ becomes a finite free coalgebra over $R'$ with comultiplication 
	\[
	\mu^{*}\colon R^{*}\to (R\otimes_{R'} R)^{*}\cong R^{*}\otimes_{R'}R^{*}.
	\]
	Extension of scalars from $R'$ to $R$ gives a map
	\[
	\Delta=\id \otimes_{R'} \mu^{*}\colon R\otimes_{R'} R^{*} \to R\otimes_{R'} R^{*} \otimes_{R'} R^{*}\cong(R\otimes_{R'} R^{*})\botimes (R\otimes_{R'} R^{*}).
	\]
	One can check that this left $R$-linear map indeed ends up in the Takeuchi product
	$(R\otimes_{R'} R^{*})\btimesr (R\otimes_{R'} R^{*})$. Translating everything back through $\End_{R'}(R)\cong R\otimes_{R'} R^{*}$ one obtains the bialgebroid structure on $\End_{R'}(R)$. The corresponding representation $\rho_{\varepsilon}$ is the inclusion $\End_{R'}(R)\hookrightarrow \End_{k}(R)$. As discussed in \cite{sweedler74} $\End_{R'}(R)$ is a Hopf algebroid over $R$.
\end{beis}
The following example is a special case of \cite[\S7, Example $A\#H$]{sweedler74} and will be important in the final section.
\begin{beis} \label{example twisted group algebras}
	Let $G$ be a group and $R$ a commutative $k$-algebra on which $G$ acts by $k$-algebra automorphisms. Consider the twisted group algebra  $R\star G\coloneqq R\otimes_k kG$. Using the short notation $rg\coloneqq r\otimes_k g$, the multiplication becomes $(rg)(r'g')=rg(r')gg'$, where $r,r'\in R, g,g'\in G$. The algebra $R\star G$, which is also known as skew group algebra or cross/crossed/smash product becomes a Hopf algebroid over $R$ when equipped with the comultiplication given by $\Delta(g)=g\botimes g$ and counit given by $\varepsilon(rg)=r$. The red map is then given by $\red_{R\star G}(rg)=r g\, \rotimes \, g^{-1}= g\, \rotimes \, g^{-1}r$. 
\end{beis}
\begin{rema} \label{remark antipodes}
	Note that in Examples \ref{example bialgebroid 1}, \ref{example bialgebroid 2}, \ref{example bialgebroid 5}, \ref{example twisted group algebras} an antipode induces the red map. More precisely there is an anti-isomorphism $S_H\colon H\to H^{\op}$ which is the identity on $R$, such that $\red_H=(\id_{H}\!\botimes S_H)\circ \Delta$, or in Sweedler notation $\red_H=h_{(1)}\rotimes S_H(h_{2})$, where $\Delta(h)=h_{(1)}\botimes h_{(2)}$ c.f.\ \Cref{example bialgebroid 1}. These are for \Cref{example bialgebroid 2} matrix transposition, for \Cref{example bialgebroid 5} the antiautomorphism $\del_i\mapsto -\del_i$ of the Weyl algebra and for \Cref{example twisted group algebras} $g\mapsto g^{-1}$. An antipode cannot exist for \Cref{example bialgebroids in finite free endo case}, unless $R$ is Frobenius over $R'$; see \cite[Theorem 12.4]{sweedler74}.
	We will see in \Cref{theorem no antipode} an explicit example for a Hopf algebroid without antipode.
\end{rema}
\section{Application to nil Hecke algebras} \label{section on Nil--Hecke basics}
This section is all about nil Hecke algebras for Coxeter systems. We recall the definitions surrounding Coxeter groups, for the theory see for instance \cite{humphreys1990}.  
\begin{nota}
	From now on we work with the ground field $k=\mR$. Additionally we fix $(W,S)$ a Coxeter system, i.e.\ a group $W$ together with a choice of generators $S\subset W$ such that $W$ has a presentation of the form $\shift{s\in S \mid s^2=1, (st)^{m_{st}}=1 \text{ for } s\neq t \in S}$ for $m_{st}=m_{ts}=\operatorname{ord}(st)\geq 2$ (possibly infinite, in which case no relation on $st$ is imposed). An expression is a finite sequence of letters from $S$ and will be denoted $\underline{w}=\underline{s_{1}\cdots s_{r}}$ with $s_i\in S$ for $1\leq i\leq r$. Given an expression $\underline{w}$ leaving out the underline means the corresponding element $w=s_{1}\cdots s_{r}\in W$. We denote by $|\underline{w}|\coloneqq r$ the length of the expression. We call an expression $\underline{w}$ of $w$ reduced if it has minimal length among the expressions yielding $w$. 
\end{nota}
Next, we recall the definitions surrounding the geometric representation $\h^{\star}$ of $(W,S)$, see for instance \cite[Theorem 1.3.11]{kumar02}.
\begin{nota} \label{definition coxeter group setup}
	Let $\h^{*}\coloneqq\bigoplus_{s\in W}\mR\!{\alpha_s}$ and call the elements $\alpha_s$ simple roots. Consider the symmetric bilinear-form on $\h^{*}$ defined by $(\alpha_s, \alpha_t)=-\cos(\frac{\pi}{m_{st}})$ (reading $m_{ss}=1$). We define for $\lambda\in \h^{*}, s\in S$ the scalar $\shift{\lambda,\alpha_s^{\vee}}\coloneqq\frac{2(\lambda,\alpha_s)}{(\alpha_s,\alpha_s)}$. For $s\in S$ and $\lambda\in \h^{\star}$ the formulas $s(\lambda)=\lambda-\shift{\lambda,\alpha_s^{\vee}}\alpha_s$ define a faithful representation of $W$ on $\h^{*}$, which we call the geometric representation of $(W,S)$. The elements of $\Phi\coloneqq \{ w(\alpha_s) \mid s\in S, w\in W\}\subset \h^{*}\setminus\{0\}$ are called roots. Let $\Phi^{+}\subseteq \Phi$ be the set of positive roots, i.e.\ those roots which are positive linear combinations of simple roots. Let $R=\Sym(\h^{*})$. We view $R$ as a graded algebra with $\h^{*}$ concentrated in degree $2$. We denote by $R^W\subseteq R$ the algebra of $W$-invariants with respect to the action of $W$ on $R$ induced by the action of $W$ on $\h^{*}$. 
\end{nota}
\begin{beis} \label{example geometric representation for Sn}
	Let $(W,S)=(S_n, \{s_i=(i,i+1) \mid 1\leq i \leq n-1\})$ and consider the action of $S_n$ on $V\coloneqq \h^{*}(\gl_n)\coloneqq \bigoplus_{i=1}^{n} \mR\!x_i$ permuting the $n$ basis vectors $x_1,\ldots, x_n$. One can identify the geometric representation $\h^{*}$ of $S_n$ with the quotient $\h^{*}(\sln)\coloneqq \h^{*}(\gl_n)/U$ by the $1$-dimensional subrepresentation $U$ spanned by $\sum_{i=1}^{n}x_i$. Under this identification the simple root $\alpha_i\coloneqq \alpha_{s_i}$ associated to $s_i=(i,i+1)$ becomes $(x_i-x_{i+1})+U$. There is no harm in replacing the geometric representation $\h^{*}$ by $\h^{*}(\gl_n)$ setting $\alpha_{i}=x_i-x_{i+1}$. In this case one defines the symmetric bilinear form $(-,-)$ in such a way that $\{x_i \mid 1\leq i\leq n\}$ becomes an orthonormal basis, in particular one has $(\alpha_i,\alpha_i)=2$. Using this choice of bilinearform, the coefficients $\shift{\alpha_i, \alpha_j^{\vee}}=\frac{2(\alpha_i, \alpha_{j})}{(\alpha_{j}, \alpha_j)}$ remain the same as for the geometric representation, moreover this redefinition of $\shift{-, \alpha_i^{\vee}}$ induces the permutation action of $S_n$ on $V$. Since we from now on never mention $(-,-)$ again and just use $\shift{-, \alpha_i^{\vee}}$, for our purposes $\h^{\star}=\h^{*}(\sln)$ and $\h^{*}(\gl_n)$ behave like they are the same. Hence, we will work with $\h^{*}(\gl_n)$ when we discuss the $S_n$ case for the nil Hecke algebra in \Cref{example nil Hecke for Sn}.
\end{beis}
\begin{rema}
	If $(W,S)$ happens to be the Weyl group of a Kac--Moody Lie algebra over $\C$ associated with a generalized Cartan matrix of size $|S|\times |S|$, there is no harm in replacing $k=\mR$ by $k=\C$ and the geometric representation $\h^{\star}$ by the dual of the Cartan, which is $(|S|+\operatorname{corank}(A))$-dimensional complex vector space. 
\end{rema}
\begin{defi} \label{definition nil Hecke}
	The nil Hecke algebra $\nh$ is the algebra generated by $\h^{*}$ (weights) and $\{\del_s\mid s\in S\}$ (nil Coxeter generators) subject to the relations
	\begin{enumerate}
		\item \label{relation comm} $\lambda \mu=\mu \lambda$ for $\lambda, \mu\in \h^{*}$. 
		\item \label{relation twisted weyl} The twisted Weyl relations
		\[
		\del_{s}\lambda=s(\lambda)\del_s + \shift{\lambda, \alpha_s^{\vee}}
		\text{ for } s\in S, \lambda\in \h^{*}.
		\] 
		\item \label{relation nil coxeter} The nil Coxeter relations 
		\[
		\del_s^2=0, \quad \underbrace{\del_s\del_t\del_s \cdots}_{m_{st}-\text{factors}}=\underbrace{\del_t\del_s\del_t \cdots}_{m_{st}-\text{factors}} \text{for } s,t\in S
		\] 
	\end{enumerate} 
	The subalgebra of $\nh$ generated by $\{\del_s\mid s\in S\}$ is called the nil Coxeter algebra and denoted by $\nc=\nc(W,S)$.
\end{defi}
\begin{beis} \label{example nil Hecke for Sn}
	Consider again the case $W=S_n$ and replace the geometric representation $\h^{\star}$ by $\h^{*}(\gl_n)$ as discussed in \Cref{example geometric representation for Sn}. In particular we set $R=\Sym(\h^{*}(\gl_n))$, which is the polynomial ring $k[x_1,\ldots, x_n]$. Writing $\del_{i}=\del_{s_i}$ this definition yields the usual type $A$ (affine) nil Hecke algebra from the literature with generators $x_1,\ldots, x_n, \del_1, \ldots, \del_{n-1}$ and relations
	\begin{gather*}
		\del_i^2=0,  \quad \del_i\del_{i+1}\del_i=\del_{i+1}\del_{i}\del_{i+1} ,
		\quad \del_i \del_j=\del_j \del_i \quad(\text{for }|i-j|\geq 2), \\
		x_ix_j=x_jx_i\quad(\text{for all } i,j), \quad \del_i x_i =x_{i+1} \del_i + 1, x_i\del_i=\del_i x_{i+1}+1, \\
		\del_i x_j =x_j \del_i \quad (\text{for } j\notin \{i,i+1\}),
	\end{gather*}
	for all indices $i,j$, such the above expressions make sense. In the literature, the nil Hecke algebra is viewed as a diagram algebra, where
	\[
	x_i=\cbox{
		\begin{tikzpicture}[tldiagram]
			\draw \tlcoord{0}{0} \lineup; 
			\draw \tlcoord{0}{2} \lineup; 
			\draw \tlcoord{0}{3} \dlineup; 
			\draw \tlcoord{0}{4} \lineup; 
			\draw \tlcoord{0}{6} \lineup; 
			\makecdots{0}{1};
			\makecdots{0}{5};
			\node at \tlcoord{0}{0} [anchor = north] {$\scriptstyle 1\vphantom{k}$};
			\node at \tlcoord{0}{2} [anchor = north] {$\scriptstyle {i-1}$};
			\node at \tlcoord{0}{3} [anchor = north] {$\scriptstyle i$};
			\node at \tlcoord{0}{4} [anchor = north] {$\scriptstyle {i+1}$};
			\node at \tlcoord{0}{6} [anchor = north] {$\scriptstyle n\vphantom{k}$};
		\end{tikzpicture}
	}, \, \del_i = \cbox{
		\begin{tikzpicture}[tldiagram]
			\draw \tlcoord{0}{0} \lineup; 
			\draw \tlcoord{0}{2} \lineup; 
			\draw \tlcoord{0}{3} \linewave{1}{1}; 
			\draw \tlcoord{0}{4} \linewave{1}{-1}; 
			\draw \tlcoord{0}{5} \lineup; 
			\draw \tlcoord{0}{7} \lineup; 
			\makecdots{0}{1};
			\makecdots{0}{6};
			\node at \tlcoord{0}{0} [anchor = north] {$\scriptstyle 1\vphantom{k}$};
			\node at \tlcoord{0}{2} [anchor = north] {$\scriptstyle {i-1}$};
			\node at \tlcoord{0}{3} [anchor = north] {$\scriptstyle i$};
			\node at \tlcoord{0}{4} [anchor = north] {$\scriptstyle {i+1}$};
			\node at \tlcoord{0}{5} [anchor = north] {$\scriptstyle {i+2}$};
			\node at \tlcoord{0}{7} [anchor = north] {$\scriptstyle n\vphantom{k}$};
		\end{tikzpicture}
	}.
	\]	
\end{beis}
Next, we gather some facts on the nil Hecke algebra in the following lemma and proposition. 
\begin{nota}
	Let $w\in W$ and let $\underline{w}=\underline{s_1\cdots s_r}$ be a expression of $w$. We write $\del_{\underline{w}}\coloneqq\del_{s_1}\cdots \del_{s_r}\in \nh$ for the nil Coxeter monomial dependent on $\underline{w}$, where by convention $\del_{\underline{\emptyset}}=1$.
\end{nota}
\begin{lemm} \label{lemma reduced expressions nil coxeter monomials}
	Let $w\in W$, let $\underline{w}=\underline{s_1\cdots s_r}$ a expression. The following are equivalent:
	\begin{enumerate}
		\item $\del_{\underline{w}}\neq 0$,
		\item $\underline{w}$ is a reduced expression.
	\end{enumerate}
	Moreover given two reduced expressions $\underline{w},\underline{w'}$ we have $\del_{\underline{w}}=\del_{\underline{w'}}$ if and only if $w=w'$.
\end{lemm}
\begin{proof}
	See \cite[Theorem 11.1.2 b), d)]{kumar02}.
\end{proof}
\begin{nota}
	By \Cref{lemma reduced expressions nil coxeter monomials}, we can define for $w\in W$ the nil Coxeter monomial $\del_{w}\coloneqq \del_{\underline{w}}$, where $\underline{w}$ is any reduced expression of $w$. The definition does not depend on this choice. 
\end{nota}
\begin{prop} \label{proposition polynomial representation nil hecke} \label{lemm the group algebra in nil Hecke}
	The following statements about $\nh$ hold:
	\begin{enumerate}
		\item \label{embeddings of twisted group algebras} Denote by $Q=\Quot(R)$ the field of fractions of $R=\Sym(\h^{*})$. The action of $W$ on $R$ extends to an action on $Q$. We have $R$-linear embeddings of $k$-algebras
		\begin{align*}
			&R\star W\hookrightarrow \nh \hookrightarrow Q\star W, \\ &s\mapsto 1-\alpha_s\del_s, \quad \del_s\mapsto \frac{1}{\alpha_s}(1-s).
		\end{align*}
		The set $\{\del_w \mid w\in W\}$ forms an $R$-basis of $\nh$ and a $Q$-basis of $Q\star W$.
		\item The assignments on generators $\nh\to \End_k(R)$
		\[
		\lambda\mapsto \rho_{\nh}(\lambda)=\lambda \cdot \blank, \quad \del_s\mapsto \odel_s\coloneqq \rho_{\nh}(\del_s)\colon f\mapsto \frac{f-s(f)}{\alpha_s}.
		\]
		define a faithful $R^{W}$-linear representation $\rho_{\nh}\colon \nh \to \End_{R^{W}}(R)$ of $\nh$ on $R$.
		The operators $\odel_s$ are called Demazure, BGG--Demazure, or simply divided difference operators, and satisfy the twisted Leibniz rule
		\begin{equation} \label{Leibniz rule}
			\odel_s(f_1)f_2+s(f_1)\odel_s(f_2)=\odel_s(f_1f_2)= \odel_s(f_1)s(f_2)+f_1\odel_s(f_2)
		\end{equation}
		for all $f_1,f_2\in R$.
		\item If $W$ is finite, one has an isomorphism of $\nh$-modules
		\[
		R\cong \nh\!e_{\triv}, \quad \text{where } e_{\triv}=\frac{1}{|W|}\sum_{w\in W}w=\frac{1}{|W|}\del_{\wo} \!\!\!\prod_{\alpha\in \Phi^{+}}\!\!\!\alpha,
		\]
		where $\wo\in W$ denotes the longest element in $W$. 
		\item If $W$ is finite, $\rho_{\nh}\colon \nh \to \End_{R^{W}}(R)$ is an isomorphism.
	\end{enumerate}
\end{prop}
\begin{proof}
	See the entirety of \cite[\S XI.1]{kumar02}. One has to carefully check that the methods of their proof, which are done for $\h$ being the Cartan subalgebra of a Kac--Moody Lie algebra, and which rely on real roots, apply in the same way to the geometric representation and its roots. 
\end{proof}
Recall that $R\star W$ is a Hopf algebroid over $R$, and $Q\star W$ is a Hopf algebroid over $Q$ with the structure discussed in \Cref{example twisted group algebras}. 
\begin{theo}[Nil Hecke bialgebroid] \label{theorem nil hecke is hopf algebroid}
	The assignments $\nh\to \nh\botimes \nh$
	\[
	\Delta_{\nh}(\lambda)\coloneqq \lambda \cdot 1 \botimes 1, \quad
	\Delta_{\nh}(\del_s)\coloneqq \del_s \botimes s + 1\botimes \del_s=\del_{s}\botimes 1 + s\botimes \del_s,
	\]
	 where $\lambda\in \h^{*}, s\in S$, extend to an algebra morphism $\Delta_{\nh}\colon \nh\to \nh \btimesr \nh$. Together with the map $\varepsilon\colon \nh\to R, \varepsilon(h)\coloneqq\rho_{\nh}(h)(1)$
	this turns $\nh$ into a cocommmutative Hopf algebroid. The map $\red_{\nh}\colon \nh \rtimesr \nh^{\op}$ is given by
	\[
	\red_{\nh}(\lambda)=\lambda 1\rotimes1, \quad \red_{\nh}(\del_s)=\del_s \rotimes s + 1\rotimes \del_s =\del_s\rotimes 1- s\rotimes \del_s.
	\] 
	The subalgebra $R\star W$ is a Hopf subalgebroid of $\nh$.
\end{theo}
To prove the theorem, we need \Cref{theorem mixed relations} below, which is about mixed versions of dihedral braid relations in $\nh$. 
\begin{nota}
	Fix $s,t\in S$ two simple reflections with $m=m_{st}<\infty$. Denote by $W_{s,t}\subseteq W$ the parabolic subgroup generated by $s,t$, which is a dihedral group of order $2m$. We fix the \cite{elias16} shorthand notation for expressions, which alternate $i$ times between $s$ and $t$
	\[
	\esi\coloneqq\underbrace{\underline{sts\cdots}}_{i \text{ letters}}, \quad \eti\coloneqq\underbrace{\underline{tst\cdots}}_{i \text{ letters}}.
	\]
	Denote by $w_{o,s,t}\in W_{s,t}$ the longest element of the parabolic subgroup $W_{s,t}$ that is the element with the two reduced expressions $\esm$ and $\etm$. For two expressions $\underline{y}$ and $\underline{w}$ we write $\underline{y}\subseteq \underline{w}$, if $\underline{y}$ is a subexpression, i.e.\ $\underline{y}$ is obtained by leaving out some letters of $\underline{w}$. We write $\iota\colon \underline{y}\hookrightarrow \underline{w}$ for an embedding, which is a fixed choice of left-out letters. An embedding $\iota$ is equivalent to a function $\chi_{\iota}\colon \{1,\ldots, |\underline{w}|\}\to\{0,1\}$, which picks out which letters are in the subexpression ($\chi_{\iota}(j)=1$) and which are not ($\chi_{\iota}(j)=0$). An expression $\underline{w}$ has exactly $2^{|\underline{w}|}$ embedded subexpressions.  
\end{nota}
\begin{beis}
	We have $\prescript{}{s}{\underline{\hat{3}}}=\underline{sts}$, $\prescript{}{s}{\underline{\hat{4}}}=\underline{stst}$ and $\prescript{}{t}{\underline{\hat{7}}}=\underline{tststst}$. If $m_{st}=5$ we have $ststs=tstst$, which implies $ \prescript{}{t}{\hat{7}}=tststst=ststsst=sts=\prescript{}{s}{\hat{3}}$. The expression $\prescript{}{s}{\underline{\hat{3}}}$ has $7$ subexpressions $\underline{\emptyset}, \underline{s},\underline{t},\underline{ts},\underline{ss},\underline{st},\underline{sts}$ and $8$ embedded subexpressions $\underline{\textcolor{red}{sts}},\underline{\textcolor{red}{st}s},\underline{s\textcolor{red}{ts}}, \underline{\textcolor{red}{s}t\textcolor{red}{s}}, \underline{\textcolor{red}{s}ts},\underline{s\textcolor{red}{t}s},\underline{st\textcolor{red}{s}}, \underline{sts}$, the red letters are those we left out of the expression. Every subexpression except $\underline{ss}$ is reduced. 
\end{beis}
\begin{defi} \label{definition mixed expressions}
	Let $\iota\colon \underline{w}\hookrightarrow \underline{w'}=\underline{s_1\cdots s_r}$ be an embedding of a subexpression. We define the mixed monomial $\mix_{\iota}\in \nh$ as the product
	\[
	\mix_{\iota}\coloneqq\mix_{\iota,1}\mix_{\iota,2}\cdots\mix_{\iota, r}\in \nh, \, \text{where }\mix_{\iota, j}=\begin{cases} 
		s_j, & \text{if } \chi_{\iota}(j)=1, \\
		\del_{s_j}, & \text{if } \chi_{\iota}(j)=0.
	\end{cases}
	\]
\end{defi}
\begin{beis} \label{example mixed expressions}
	Consider $\prescript{}{s}{\hat{\underline{7}}}=\underline{stststs}$. For  $\iota_1\colon\underline{\textcolor{red}{s}t\textcolor{red}{st}sts}\hookrightarrow \prescript{}{s}{\hat{\underline{7}}}$ we have $\mix_{\iota_1}=\del_s t \del_s \del_t sts$. For $\iota_2\colon \underline{\textcolor{red}{s}t\textcolor{red}{s}ts\textcolor{red}{ts}}\hookrightarrow \prescript{}{s}{\hat{\underline{7}}}$ we have $\mix_{\iota_2}=\del_s t \del_s ts \del_t\del_s$.
\end{beis}
\begin{theo}[Mixed dihedral braid relations] \label{theorem mixed relations}
	Let $s,t\in S$ two simple reflections with $m=m_{st}<\infty$. Let $w\in W_{s,t}, w\neq w_{o,s,t}$ be a fixed  element. Then the relation
	\[
	\rR_w\colon \sum_{\substack{\iota\colon \underline{w}\subseteq \esm \\ \underline{w} \text{ is reduced}}}\mix_{\iota}=\sum_{\substack{\iota\colon \underline{w}\subseteq \etm \\ \underline{w} \text{ is reduced}}}\mix_{\iota}
	\]
	holds in $\nh$.
\end{theo}
\begin{beis} \label{example new relations}
	The relation $\rR_{1}$ is always $\del_{\esm}=\del_{\etm}$.
	The other relations are mixed versions of the braid relation in the dihedral group $W_{s,t}$. For instance:
	\begin{enumerate}
		\item For $m_{st}=2$ we obtain the relations
		\[
		\rR_{s}\colon s\del_t=\del_ts, \quad \rR_{t}\colon \del_s t=t \del_s
		\]
		\item For $m_{st}=3$ we obtain the type $A$ relations
		\[
		\rR_{s}\colon s\del_t\del_s+\del_s\del_ts=\del_ts\del_t, \quad \rR_{st}\colon st\del_s=\del_tst.
		\]
		as well as $\rR_{t}$ and $\rR_{ts}$, which are the same as the above relations with $s,t$ replaced. 
		\item For $m_{st}=4$ we obtains the type $B/C$ relations
		\begin{gather*}
			\rR_{s}\colon s\del_t\del_s\del_t+\del_s\del_ts\del_t=\del_ts\del_t\del_s+\del_t\del_s\del_ts, \\ \rR_{st}\colon st\del_s\del_t+s\del_t\del_s t+\del_s \del_t st=\del_t st\del_s, \quad \rR_{sts}\colon sts\del_t=\del_t sts,
		\end{gather*}
		and as before, the ones obtained by swapping $s$ and $t$.
	\end{enumerate}
\end{beis}
\begin{proof}[Proof of \Cref{theorem mixed relations} on mixed relations]
	Our proof relies on the embedding $\nh\subseteq Q\star W$ and on the $Q$-linear comultiplication on $Q\star W$. We write $\cotimes\coloneqq\botimes_Q$ to distinguish it from $\botimes=\botimes_R$. First, we compute that the comultiplication 
	$\Delta\colon Q\star W\to (Q\star W) \cotimes (Q\star W)$ maps $\del_s$ to $\del_s \cotimes s + 1\cotimes \del_s$. Indeed one has
	\begin{align*}
		\del_s \cotimes s + 1\cotimes \del_s&= \frac{1}{\alpha_s}(1-s)\cotimes s + 1\cotimes \frac{1}{\alpha_s}(1-s) \\
		&=\frac{1}{\alpha_s} (1\cotimes s- s\cotimes s + 1\cotimes 1 - 1\cotimes s)
		&=\frac{1}{\alpha_s}(1\cotimes 1-s\cotimes s)=\Delta(\del_s).
	\end{align*}
	Let $s,t\in S$ with $m\coloneqq m_{st}< \infty$. In $(Q\star W)\cotimes (Q\star W)$ we have the equality
	\[
	\Delta(\underbrace{\del_s\del_t\cdots}_{m_{st} \text{ factors}})=\Delta(\underbrace{\del_t\del_s\cdots}_{m_{st} \text{ factors}}).
	\]
	Using that $\Delta$ is an algebra morphism into $(Q\star W)\cotimes(Q\star W)$ this implies that
	\[
	\underbrace{\Delta(\del_s)\Delta(\del_t)\cdots}_{m_{st} \text{ factors}} =\underbrace{\Delta(\del_t)\Delta(\del_s)\cdots}_{m_{st} \text{ factors}}
	\]
	Using the formula for $\Delta(\del_s)$ and $\Delta(\del_t)$ and multiplying everything out on both sides, one obtains
	\[
	\sum_{\iota\colon \underline{w}\subseteq \esm}\del_{\underline{w}}\cotimes \mix_{\iota}=\sum_{\iota\colon \underline{w}\subseteq \etm}\del_{\underline{w}}\cotimes \mix_{\iota},
	\]
	where the sums runs over all embedded subexpressions $\underline{w}$ of $\esm$ respectively $\etm$.
	By \Cref{lemma reduced expressions nil coxeter monomials}, we have $\del_{\underline{w}}=0$, unless $\underline{w}$ is a reduced expression in $W$. Hence the equality becomes 
	\[
	\del_{\esm}\cotimes \!\sm+\sum_{w\in W_{s,t}}\del_{w}\,\cotimes \!\!\!\!\!\!\sum_{\substack{\iota\colon \underline{w}\subseteq \esm \\ \underline{w} \text{ is reduced}}}\mix_{\iota}= \del_{\etm}\cotimes \!\tm +\sum_{w\in W_{s,t}}\del_{w}\,\cotimes\!\!\!\!\!\! \sum_{\substack{\iota\colon \underline{w}\subseteq \etm \\ \underline{w} \text{ is reduced}}}\mix_{\iota}
	\]
	after reordering summands.
	The first summands agree because $\del_{\esm}=\del_{\etm}$ and $\sm=\tm$ in $Q\star W$. Since $\{\del_w \mid w\in W\}$ is a $Q$-basis of $Q\star W$ by \Cref{proposition polynomial representation nil hecke} we can compare coefficients to deduce that for $w\neq w_{o,s,t}$ 
	\[
	\sum_{\substack{\iota\colon \underline{w}\subseteq \esm \\ \underline{w} \text{ is reduced}}}\mix_{\iota}=\sum_{\substack{\iota\colon \underline{w}\subseteq \etm \\ \underline{w} \text{ is reduced}}}\mix_{\iota}
	\]
	in $Q\star W$ and in particular in $\nh$.
\end{proof}
Now we can prove \Cref{theorem nil hecke is hopf algebroid}, i.e. that nil Hecke algebras form Hopf algebroids.
\begin{proof}[Proof of \Cref{theorem nil hecke is hopf algebroid}]
	We first check $\del_s \botimes s + 1\botimes \del_s=\del_{s}\botimes 1 + s\botimes \del_s$. By $R$-linearity of $\botimes$ we have
	\begin{align*}
		\del_s\botimes 1 + s \botimes \del_s &= \del_s \botimes 1 + (1-\alpha_s\del_s)\botimes \del_s \\
		&=\del_s \botimes 1 + 1\botimes \del_s - \del_s\botimes \alpha_s\del_s = \del_s \botimes s + 1\botimes \del_s.
	\end{align*}
	Since $R$ is generated by $\h^*$, it is enough to check the Takeuchi condition with elements of $\lambda\in \h^*$. We calculate
	\begin{align*}
		\del_{s}\lambda\botimes 1 + s \lambda \botimes \del_{s} &= (s(\lambda)\del_s+\shift{\lambda, \alpha_s^{\vee}})\botimes 1 + s(\lambda)s\botimes \del_s \\
		&= \del_s\botimes s(\lambda)+1\botimes \shift{\lambda, \alpha^{\vee}} + s\botimes s(\lambda) \del_s.
	\end{align*}
	Now we use the twisted Weyl relation on the rightmost summand and obtain
	\begin{align*}
		&\hphantom{=}\del_s\botimes s(\lambda)+1\botimes \shift{\lambda, \alpha_s^{\vee}}+ s\botimes(\del_s\lambda-\shift{\lambda, \alpha_s^{\vee}}) \\
		&=\del_s\botimes (\lambda-\shift{\lambda, \alpha_s^{\vee}}\alpha_s)+1\botimes \shift{\lambda, \alpha_s^{\vee}} + s\botimes(\del_s\lambda-\shift{\lambda, \alpha_s^{\vee}}) \\
		&=\del_s\botimes \lambda + s\botimes \del_s\lambda,
	\end{align*}
	where we used in the last step that 
	\[
	-\del_s\botimes \shift{\lambda, \alpha_s^{\vee}}\alpha_s = - \alpha_s \del_s\botimes \shift{\lambda, \alpha_s^{\vee}}=(s-1)\botimes \shift{\lambda, \alpha_s^{\vee}}. \]
	When checking that the comultiplication $\Delta_{\nh}\colon \nh \to \nh \btimesr \nh$ is well-defined, the only non-trivial part is 
	\[
	\underbrace{\Delta_{\nh}(\del_s)\Delta_{\nh}(\del_t)\cdots}_{m_{st} \text{ factors}} =\underbrace{\Delta_{\nh}(\del_t)\Delta_{\nh}(\del_s)\cdots}_{m_{st} \text{ factors}},
	\]
	which is equivalent by the same computation as in the proof of \Cref{theorem mixed relations} to the following equality in $\nh \botimes \nh$
	\[
	\del_{\esm}\botimes \!\sm+\sum_{w\in W_{s,t}}\del_{w}\,\botimes \!\!\!\!\!\!\sum_{\substack{\iota\colon \underline{w}\subseteq \esm \\ \underline{w} \text{ is reduced}}}\mix_{\iota}= \del_{\etm}\botimes \!\tm +\sum_{w\in W_{s,t}}\del_{w}\,\botimes\!\!\!\!\!\! \sum_{\substack{\iota\colon \underline{w}\subseteq \etm \\ \underline{w} \text{ is reduced}}}\mix_{\iota}.
	\]
	By \Cref{theorem mixed relations}, this is true. The rest of the definition is trivial, note that $\varepsilon(\del_s)=0$ for all $s\in S$ to check that $\nh, \Delta, \varepsilon$ is a coalgebra over $R$. To check that $\nh$ is a Hopf algebroid, one checks that the assignments we wrote extend to a well-defined algebra morphism $\red_{\nh}\colon \nh \rtimesr \nh^{\op}$, which induces $\Gal_{\nh}^{-1}$.
\end{proof}
\begin{rema}
	The descent theorem \cite[\S5.5]{kraehmermahaman24} gives a general condition on how to restrict a Hopf algebroid structure on some algebra to a given subalgebra. In our setting, the `original' Hopf algebroid is $Q\star W$, and the subalgebra is the nil Hecke algebra. From this point of view, the freeness argument behind the proof of the mixed relations \Cref{theorem mixed relations} is there to show the conditions written in \cite[Lemma 5.13]{kraehmermahaman24}. Indeed, \Cref{theorem nil hecke is hopf algebroid} can be seen as a special case of the descent theorem.
\end{rema}
\begin{koro} \label{corollary equivalence is monoidal}
	Assume $W$ is finite. Then the symmetric monoidal structure on $\nh\Mod$ coming from the cocommutative $R$-bialgebroid structure on $\nh$ corresponds under the equivalence
	\[
	R\otimes_{R^{W}}\colon R^{W}\Mod \to \nh\Mod
	\] 
	to the symmetric monoidal product $\otimes_{R^{W}}$ on $R^{W}\Mod$.
\end{koro}
\begin{proof}
	An easy and long calculation shows that the comultiplication $\Delta$ of $\nh$ corresponds under the isomorphism $\rho_{\nh}\colon \nh \to \End_{R^{W}}(R)$ to the bialgebroid structure from \Cref{example bialgebroids in finite free endo case}, which is the one corresponding to the Morita equivalence. 
\end{proof}

\begin{koro} \label{theorem no antipode}
	Nil Hecke algebras are Hopf algebroids without an antipode. More precisely $\red_{\nh}\neq (\id_{\nh}\!\botimes S)\circ \Delta$ for any $k$-algebra morphism $S\colon \nh \to \nh^{\op}$, which is the identity on $R$.
\end{koro}
\begin{proof}
	To simplify notation we just consider the case where $W=S_2=\{1,s\}$ and $\h^{*}=\h^{*}(\sl2)=\mR\! \alpha$. We write $\alpha\coloneqq\alpha_s$ and $\del=\del_s$ for the generators of $\nh$, in particular $s=1-\alpha\del=\del\alpha-1$. The proof has three steps. First we calculate that $\red_{\nh}(s)=s\rotimes s$. Assuming that an antipode $S$ exists, we show in the second step that $S(s)=s$. Finally, we let $S(\del_s)\alpha$ act on $R$ and use degree reasons to get a contradiction. In step 1, we calculate
	\begin{align*}
		\red_{\nh}(1-\alpha\del)&=1\rotimes 1 - \alpha (\del \rotimes s + 1\rotimes \del) \\
		&=1\rotimes(1+s-s) -\alpha\del \rotimes s + \alpha 1\rotimes \del \\
		&= s\rotimes s + 1\rotimes (1-s) + \alpha 1\rotimes \del = s\rotimes s,
	\end{align*}
	where we used $1-s=\alpha\del$. Now assume an antipode $S$ inducing $\red_{\nh}$ exists. We know that $\Delta(s)=s\botimes s$, hence $s\rotimes s=\red_{\nh}(s)=s\rotimes S(s)$. In order to conclude that $s=S(s)$ recall that $\nh$ is free as a right $R$-module with basis $1,\del$ and hence $\nh \rotimes \nh^{\op}$ is a free right $\nh$-module with basis $1\rotimes 1, \del\rotimes 1$. We have $s=\del \alpha -1$, which leads to
	\[
	(\del \alpha -1 ) \rotimes s= (\del \alpha -1 ) \rotimes S(s).
	\]
	In particular we compare coefficients in $1\cdot R\, \rotimes \nh^{\op}$ to conclude that $s=S(s)$. To come to a contradiction, consider 
	$s=S(s)=S(1-\alpha \del)=1-S(\del)\alpha$, which implies $1-s=S(\del)\alpha$. Now consider the action of both sides on $\alpha\in R$. We have $
	2\alpha=(1-s)(\alpha)=(S(\del)\alpha)(\alpha)=S(\del)(\alpha^2)=\alpha^2S(\del)(1)$.
	Now $S(\del)(1)\in R$ is some polynomial in $\alpha$ which multiplied with $\alpha^2$ yields $2\alpha$, which is a contradiction.
\end{proof}
\begin{rema}
	In the $S_2$-case we have $\nh\cong\End_{k[\alpha^2]}(k[\alpha])$ which under appropiate identifications becomes \Cref{example bialgebroids in finite free endo case}. In \cite[Theorem 12.4]{sweedler74}, it is discussed that if an antipode on $\End_{R'}(R)$ exists, then $R$ has to be Frobenius as an algebra over $R'$. This is the case here: there is no isomorphism of $k[\alpha]$-modules $k[\alpha]\to k[\alpha]^{*}=\Hom_{k[\alpha^2]}(k[\alpha],k[\alpha^2])$, since $k[\alpha]^{*}$ contains torsion elements for $\alpha$ like $1^{\vee}\colon 1\mapsto 1, \alpha\mapsto 0$. However, the moral point of our proof is that $Q\star W$ has an antipode (as Hopf algebroid over $Q$), which does not map $\nh$ into itself. Indeed the antipode is the algebra morphism $Q\star W \to (Q \star W)^{\op}$, which maps maps $s$ to $s$ and $\alpha$ to $\alpha$, hence $\del=\frac{1}{\alpha}(1-s)$ is mapped to $(1-s)\frac{1}{\alpha}$, which is not in $\nh$, since it does not preserve $R=k[\alpha]\subset Q=k(\alpha)$.
\end{rema}

	\ihead{}
	\chead{} 
	\ohead{References}  
	\printbibliography
\end{document}